\documentclass{article}

\usepackage{amsthm}
\usepackage{amsmath}
\usepackage{amssymb}

\usepackage{paralist}

\newtheorem{lemma}{Lemma}

\theoremstyle{definition}
\newtheorem{algorithm}{Algorithm}

\newcommand{\N}{\mathbb{N}}
\newcommand{\R}{\mathbb{R}}

\title{Exploring Mount Neverest}
\author{Michiel de Bondt}

\begin{document}

\maketitle

In one of the columns in the series `Perplexities' in 1922, 
Henry Ernest Dudeney formulated the following problem:

\begin{quotation}
Professor Walkingholme, one of the exploring party, was allotted the
special task of making a complete circuit of the base of the mountain at
a certain level.  The circuit was exactly a hundred miles in length and
he had to do it all alone on foot.  He could walk twenty miles a day,
but he could only carry rations for two days at a time, the rations for
each day being packed in sealed boxes for convenience in dumping.  He
walked his full twenty miles every day and consumed one day's ration as
he walked.  What is the shortest time in which he could complete the
circuit?
\end{quotation}

This problem can be found in the book `536 Puzzles \& Curious Problems' from
Henry Ernest Dudeney, edited by Martin Gardner. But I did not find
an optimal solution to the problem in the literature. Albeit Martin Gardner is making fun
on it, the right interpretation of the problem is not clear to me at all.
Let us first formulate some ways to tackle the problem.

\section{Solution of the problem}

One way to make the circuit is doing the same as with a straight distance of $100$ miles. Doing that and minimizing the traveled distance, one can prove that $82\frac{8}{55} \cdot 20$ miles need to be traveled. So for this approach, $82\frac{8}{55}$ is a lower bound of the required amount of time in days. It seems a better approach to walk two round trips to the $50$ miles distance point on the other side of the mountain, since one can prove that $41\frac{1}{4} \cdot 20$ miles need to be traveled for that. So for this approach, $41\frac{1}{4}$ is a lower bound of the required amount of time in days. But Dudeney found the following solution:

\begin{algorithm} \label{dudeney} ($23\frac12$ days)
\begin{enumerate}
\item Dump 5 rations at 90-mile point and return to base (5 days).
\item Dump 1 at 85 and return to 90 (1 day).
\item Dump 1 at 80 and return to 90 (1 day).
\item Dump 1 at 80, return to 85, pick up 1 and dump at 80 (1 day).
\item Dump 1 at 70 and return to 80 (1 day).
\item Return to base (1 day). We have thus left one ration at 70 and one at
      90.
\item Dump 1 at 5 and return to base (1 day).  If he must walk 20 miles
      he can do so by going to 10 and returning to base.
\item Dump 4 at 10 and return to base (4 days).
\item Dump 1 at 10 and return to 5; pick up 1 and dump at 10 (1 day).
\item Dump 2 at 20 and return to 10 (2 days).
\item Dump 1 at 25 and return to 20 (1 day).
\item Dump 1 at 30, return to 25, pick up 1 and dump at 30 (1 day).
\item March to 70 (2 days).
\item March to base ($1 \frac12$ days).
\end{enumerate}
\end{algorithm}

In step 7, only 10 miles are walked effectively, so 10 miles are wasted. A way to avoid this is the following. Professor Walkinghome starts walking in the middle of the day, so half of the box which is used that day has already been consumed. With the other half, Professor Walkinghome performs step 7 first. After that, the other steps follow. Then Professor Walkinghome starts and ends in the middle of the day, and the time in between measures 23 days.

We can divide Dudeney's solution in three parts:
\begin{description}
\item[Part A:]
a round trip to the 70 miles point from the base, in which boxes are dumped
on positions 70 and 90 for later use (steps 1 to 6, $10$ days).
\item[Part B:]
a one-way trip to the 70 miles point from the base around the mountain
(steps 7 to 13, $11 \frac12$ or $12$ days).
\item[Part C:]
walking from the 70 miles point to the base as the 100 miles point, using the boxes dumped in part A (step 14, $1 \frac12$ 
days).
\end{description}
If 12 days are counted for part B, then it does not seem a good idea to waste 10 miles. So let us count $11\frac12$ days for part B for now. Only for part C, it is immediately clear that is optimal. Later, we will see that part A is optimal as well, but part B can be improved, to $11 \frac17$ days exactly.

Before trying to find a better solution, it is always a good idea what others did. When I searched the internet with google, I found the homepage of a youngster called Nightvid Cole, who presents better solutions than that of Dudeney. In the first one, he presents a solution in which $22 \frac{10}{11} \cdot 20$ miles are traveled. The solution can easily be adapted to a solution in which Professor Walkinghome starts on a specific moment during the day and walks $22 \frac{10}{11}$ days. Although Cole's solution is better than that of 
Dudeney, only part C of it is optimal.

The improvement of Cole is that he dumps boxes on positions $70\frac{10}{11}$ and $90 \frac{10}{11}$ rather than $70$ and $90$. Now part B must reach further, whence it takes $12$ days now. But probably, Cole reasoned the other way around: he reserved 12 days for part B and then thought out an according solution. That was a very good idea, but since he failed to optimize parts A and B, his solution is not optimal. So I optimized parts A and B, which resulted in the following solution, a solution that turned out to be optimal later on.

\begin{algorithm} \label{best} ($22 \frac 9{16}$ days)
\begin{enumerate} 
\item Dump one ration at $98 \frac34$ point and return to base ($\frac18$ day).
\item Dump one ration at $97 \frac12$, return to $98 \frac34$, pick up one,
      dump at $91 \frac14$ and return to base (1 day).
\item Dump one ration at $93 \frac34$, return to $97 \frac12$, pick up one,
      dump at $93 \frac34$ also and return to base (1 day).
\item Dump two rations at $90$ and return to base (2 days).
\item Dump one ration at $86 \frac78$ and return to $93 \frac34$ (1 day).
\item Dump one ration at $82 \frac12$, return to $86 \frac78$, pick up one,
      dump at $86 \frac14$ and return to $90$ (1 day).
\item Dump one ration at $80$, return to $86 \frac14$, pick up one and get to
      $82 \frac12$ (1 day).
\item Dump one ration at $71 \frac14$ and return to $80$ (1 day).
\item Return to base (1 day). We have thus left one ration at $71 \frac 14$ 
      and one at $91 \frac14$.
\item Dump five rations at $10$ and return to base (5 days).
\item Dump one ration at $12 \frac12$, return to $10$, pick up one,
      dump at $12 \frac12$ also and return to $10$ (1 day).
\item Dump one ration at $20$ and return to $10$ (1 day).
\item Dump one ration at $20 \frac58$, return to $20$, pick up one,
      dump at $20 \frac58$ also and return to $12 \frac12$ (1 day).      
\item Dump one ration at $26 \frac9{16}$ and return to $20 \frac58$ (1 day).
\item Dump one ration at $31 \frac14$, return to $26 \frac9{16}$, pick up one
      and get to $31 \frac14$ (1 day).
\item March to $71 \frac14$ (2 days).
\item March to base ($1 \frac7{16}$ day).
\end{enumerate}
\end{algorithm}

If you look at the above algorithm, then one thing immediately strikes:
it would have been nicer if the circuit would have been 160 kilometers, with
a unit distance of 32 kilometers a day. Dudeney only considered solutions 
from which the eating and turning points were a multiple of 5 miles. This is
however impossible for a solution of $22 \frac 9{16}$ days (for $\frac{20}5
\cdot 22 \frac 9{16}$ is not integral). 

In Dudeney's solution, we counted $11\frac12$ miles for part B by allowing Professor Walkinghome to start on a specific moment during the day. But if Professor Walkinghome must start at dawn on some day, then this does not work. With the above solution, $17\frac12$ miles are wasted on the first day. One can use these miles to increase the points where boxes are dumped in part A for use in part C. Nightvid Cole found a solution of $23 \frac13$ days in this context, but again, parts A and B are not optimal. So I optimized parts A and B again, which resulted in the following solution, a solution that turned out to be optimal later on.

\begin{algorithm} \label{dawn} ($23\frac{25}{116}$ days)
\begin{enumerate}
\item Dump one ration at $8 \frac{18}{29}$ and return to base ($\frac{25}{29}$
      days).
\item Dump two rations at $99 \frac9{29}$ and return to base ($\frac4{29}$ 
      days).
\item Dump one ration at $96 \frac{26}{29}$, return to $99 \frac9{29}$,
      pick up two in turn, dump both at $95 \frac{25}{29}$ and 
      return to base (1 day).
\item Dump one ration at $90$ and return to base (1 day).
\item Dump one ration at $88 \frac{28}{29}$, return to $90$, pick up one, dump
      at $88 \frac{28}{29}$ also and return to $95 \frac{25}{29}$ (1 day).
\item Dump one ration at $82 \frac{12}{29}$ and return to $88 \frac{28}{29}$
      (1 day).
\item Dump one ration at $75 \frac{20}{29}$ and return to $82 \frac{12}{29}$
      (1 day).
\item Return to $96 \frac{26}{29}$, pick up one, dump at $95 \frac{20}{29}$
      and return to base (1 day). We have thus left one ration at 
      $8 \frac{18}{29}$, one at $75 \frac{20}{29}$ and another one at 
      $95 \frac{20}{29}$.
\item Dump one ration at $9 \frac{9}{29}$, return to $8 \frac{18}{29}$, pick 
      up one, dump at $9 \frac{9}{29}$ also and return to base (1 day).    
\item Dump five rations at $10$ and return to base (5 days).
\item Dump one ration at $12 \frac{19}{58}$, return to $10$, pick up one, 
      dump at $12 \frac{19}{58}$ also and return to $9 \frac{9}{29}$ (1 day).    
\item Dump one ration at $19 \frac{19}{29}$ and return to $10$ (1 day).
\item Dump one ration at $19 \frac{24}{29}$, return to $19 \frac{19}{29}$, 
      pick up one, dump at $19 \frac{24}{29}$ also and return to $10$ (1 day).    
\item Dump one ration at $21 \frac{19}{116}$ and return to $12 \frac{19}{58}$ 
      (1 day).
\item Dump one ration at $23 \frac{18}{29}$, return to $21 \frac{19}{116}$, 
      pick up one, dump at $23 \frac{18}{29}$ also and return to 
      $19 \frac{24}{29}$ (1 day).
\item Dump one ration at $31 \frac{21}{29}$ and return to $23 \frac{18}{29}$
      (1 day).
\item Dump one ration at $35 \frac{20}{29}$, return to $31 \frac{21}{29}$, 
      pick up one and get to $35 \frac{20}{29}$ (1 day).
\item March to $75 \frac{20}{29}$ (2 days).
\item March to base ($1 \frac{25}{116}$ days).      
\end{enumerate}      
\end{algorithm}

Notice that the above algorithm not only has ugly positions, but also that both part A and part B are partially done on the first day.

In the 2010 version of this paper, I assumed that no rations are exchanged between different parts of the exploration. For instance, in the case of two round tips, I assumed that no rations of the first round trip are used in the second round trip. But on Mathematics Stack Exchange, it has been pointed out by Forero that such assumptions undermine optimality proofs \cite{forero}. This has led to the following lemma.

\begin{lemma} \label{ABC}
We may assume that the exploration can be divided into parts A, B and C as above, such that the boundary of part A and part B is not crossed by any ration box, sealed or unsealed. Furthermore, we may assume that part C is a trip from the boundary of part A and part B straight to the starting point.
\end{lemma}

\begin{proof}
Assume that the exploration takes less than $25$ days, just as with the three given solutions. In the next section, we will prove that there is a position $z$ which is visited only one time and crossed. So Professor Walkingholme reaches position $z$ from one side, say the B-side, and leaves position $z$ on the other side, say the A-side.

Suppose that on the moment that Professor Walkingholme is on position $z$, the ration boxes which are on the A-side lie on positions $x_1, x_2, \ldots, x_l$, such that for all $i$, position $x_{i+1}$ is farther in the A-side from position $z$ than position $x_i$.

Professor Walkingholme is on position $z$ with at most $2$ ration boxes. Let $f$ be the total amount of food in rations which Professor Walkingholme carries. Then $0 < f < 2$. The distance within the A-side between position $z$ and position $x_i$ is at most $20 \cdot (i - 1 + f)$, because Professor Walkingholme can use at most $(i - 1 + f)$ rations to reach $x_i$, namely the $f$ rations on position $z$ and the ration boxes on positions $x_1, x_2, \ldots, x_{i-1}$. 

Professor Walkingholme could do the following instead. He leaves position $z$ on the A-side and walks straight to the starting point. Each time that he has no food left, he unseals a ration box, say on position $y_1, y_2, \ldots, y_{l'}$. It is indeed possible to get ration boxes on positions $y_1, y_2, \ldots, y_{l'}$ instead of positions $x_1, x_2, \ldots, x_l$, because for each $i \le l'$, position $y_i$ is at least as far on the A-side from position $z$ as position $x_i$, so position $x_i$ is at least as far on the A-side from the starting point as position $y_i$.  Professor Walkingholme does not need to reach farther than $y_1$ in part A, so $y_1$ becomes the boundary of part A and part B.
\end{proof}

\section{Exploring as a jeep}

There are basically two ways to perform the exploration. The first way is that Professor Walkingholme goes around the mountain, so the exploration has winding number $1$ or $-1$ with respect to the mountain top. The second way is that Professor Walkingholme does not go around the mountain, so the exploration has winding number $0$ with respect to the mountain top. 

\begin{lemma}
If $r$ rations cross position $x$ (counter)clockwise and the exploration has winding number $(-)w$, then Professor Walkingholme visits $x$ at least
$$
\max\{2\lceil r/2 \rceil + w, -w\}
$$
times.
\end{lemma}

\begin{proof}
Suppose that $r$ rations cross position $x$ clockwise and the exploration has winding number $w$ (the other case is similar). Suppose first that $-w \ge \lceil r/2 \rceil$. Then Professor Walkingholme must both enter and leave position $x$ clockwise $-w$ times. So position $x$ is visited at least $-w$ times.

Suppose next that $-w \le \lceil r/2 \rceil$. To get the transportations on position $x$ done, Professor Walkingholme must both enter and leave position $x$ clockwise $\lceil r/2 \rceil$ times, because Professor Walkingholme can carry at most $2$ rations. To get the winding number correct, Professor Walkingholme must both enter and leave position $x$ counterclockwise $\lceil r/2 \rceil + w$ times in addition. So position $x$ is visited at least $2\lceil r/2 \rceil + w$ times.
\end{proof}

Suppose that $r \in \R$ rations leave the starting point clockwise, where we do not count parts of rations which are not used for walking. Let $p(t)$ be the point such that Professor Walkingholme walks $20 \cdot (t+r)$ miles between the starting point and $p(t)$ on the clockwise side of the starting point. Then in some sense, $-t$ used rations cross $p(t)$ clockwise, since $t + r$ rations are used between the starting point and $t$ on the clockwise side of the starting point. But there may be used rations which cross $t$ counterclockwise as well. $t$ is the amount of used rations which crosses $p(t)$ counterclockwise minus the amount of used rations which crosses $p(t)$ clockwise.

Among other things, a theorem of Banach \cite[Th\'{e}or\`{e}me 2]{banach} says the following: if each position of an interval is visited $k$ times, then the amount of walking in that interval is at least $k$ times the length of that interval. See \cite{dozen} for more information about this. Let $\operatorname{cd}(a,b)$ denote the distance if one walks from $a$ clockwise to $b$.

Suppose first that the exploration has winding number $0$. Then for $t \ge 0$, $p(-t)$ and $p(t)$ are visited at least $2 \lceil t/2 \rceil$ times. So for $k \in \N$ and $t > 2k$, we have
\begin{align*}
\operatorname{cd}\big(p(-t),p(-2k)\big) &\le 20\cdot\frac{t-2k}{2(k+1)} &
\operatorname{cd}\big(p(2k),p(t)\big) &\le 20\cdot\frac{t-2k}{2(k+1)}
\end{align*}
In particular
\begin{align*}
\operatorname{cd}\big(p(-2(k+1)),p(-2k)\big) &\le 20\cdot\frac{2}{2k+2} &
\operatorname{cd}\big(p(2k),p(2(k+1))\big) &\le 20\cdot\frac{2}{2k+2}
\end{align*}
From
$$
20 \cdot \Big(\frac{1\frac2{5}}{14} + \frac2{12} + \frac2{10} + \frac2{8} + \frac2{6} + \frac2{4} + \frac2{2} + \frac2{2} + \frac2{4} + \frac2{6} + \frac2{8} + \frac2{10} + \frac2{12}\Big) = 100
$$
one can infer that at least $25\frac25$ days are necessary to explore Mount Neverest with winding number $0$.

Suppose next that the exploration has winding number $-1$. Then for $t \ge 0$, $p(-t)$ is visited at least $2 \lceil t/2 \rceil - 1$ times and $p(t)$ is visited at least $2 \lceil t/2 \rceil + 1$ times. So for $k \in \N$ and $t > 2k$, we have
\begin{align*}
\operatorname{cd}\big(p(-t),p(-2k)\big) &\le 20\cdot\frac{t-2k}{2(k+1)-1} & \operatorname{cd}\big(p(2k),p(t)\big) &\le 20\cdot\frac{t-2k}{2(k+1)+1}
\end{align*}
In particular
\begin{align*}
\operatorname{cd}\big(p(-2(k+1)),p(-2k)\big) &\le 20\cdot\frac{2}{2k+1} &
\operatorname{cd}\big(p(2k),p(2(k+1))\big) &\le 20\cdot\frac{2}{2k+3}
\end{align*}
From
$$
20 \cdot \Big(\frac{\frac{23}{35}}{9} + \frac2{7} + \frac2{5} + \frac2{3} + \frac2{1} + \frac2{3} + \frac2{5} + \frac2{7} + \frac2{9}\Big) = 100
$$
one can infer that at least $16\frac{23}{35}$ days are necessary to explore Mount Neverest with winding number $-1$ or $1$.

Suppose finally that there does not exist a position $z$ which is visited one time and crossed. If the winding number is $0$, then the exploration takes more than $25$ days. If the winding number is nonzero, then every position except maybe the starting point is crossed, so there does not exist a position $z$ which is visited one time. Assume without loss of generality that the winding number is $-1$. Then for $k \in \N$ and $t > 2k$, we have
\begin{align*}
\operatorname{cd}\big(p(-t),p(-2k)\big) &\le 20\cdot\frac{t-2k}{\max\{2(k+1)-1,2\}}
\end{align*}
In particular
\begin{align*}
\operatorname{cd}\big(p(-2(k+1)),p(-2k)\big) &\le 20\cdot\frac{2}{\max\{2k+1,2\}}
\end{align*}
From
$$
20 \cdot\Big( \frac1{13} + \frac2{11} + \frac2{9} + \frac2{7} + \frac2{5} + \frac2{3} + \frac2{2} + \frac2{3} + \frac2{5} + \frac2{7} + \frac2{9} + \frac2{11} + \frac2{13}\Big) = 94\frac{7858}{9009}
$$
one can infer that more than $25$ days are necessary to explore Mount Neverest if there does not exist a position $z$ which is visited one time and crossed.

Notice that the above lower bounds hold as well, if Professor Walkingholme is replaced by a jeep as in \cite{jeep} and \cite{dozen}, which can hold $2$ units of fuel in its tank, and which can ride $20$ mile per unit of fuel. For the jeep, the lower bounds of $25\frac25$ days and $16\frac{23}{35}$ days can be attained. The proof of that is left as an exercise to the reader.

We can make assumptions such that the above lower bounds can be attained for Professor Walkingholme. In the previous section, we assumed that Professor Walkingholme unseals the ration which is used that day at dawn. But we can also allow more rations to be used partially on some day. In addition we can assume one of the following two options.
\begin{compactitem}

\item
Professor walkingholme can eat enough food to walk $2$ yards, so he does not need to keep a ration box with him to eat from when he walks these $2$ yards.

\item
Food may be transferred freely between ration boxes. In addition to the ration boxes, there is a small box which can contain food for $1$ yard, which Professor Walkingholme can carry in addition to the $2$ ration boxes, but the total amount of food which is carried does not exceed $2$ rations.

\end{compactitem}
In the first option, the total amount of food which Professor Walkingholme can carry is $2$ ration boxes plus the amount of food to walk $2$ yards, which is too much. But if we subtract from each box the amount of food to walk $1$ yard, then Professor Walkingholme does not carry too much food any more, and a minimal amount of extra time will be needed.

Suppose that Professor Walkingholme is with $2k$ ration boxes on some point. Then he can carry those boxes one yard farther in $k$ turns, taking $2k-1$ yards altogether. He can manage that the food needed for that comes from only one ration box as follows. With the first option, he sees the ration box he uses every $2$ yards of walking, which is enough. With the second option, he uses the small box for the $k-1$ walks back to his first position, and uses food of any ration box for the $k$ forward walks and for filling the small box. Afterwards, he transfers food between the $2k$ ration boxes such that one ration box is effectively used.

If Professor Walkingholme can carry $2k$ boxes one yard farther as described above, then the lower bounds of $25\frac25$ days and $16\frac{23}{35}$ days can indeed be attained, using $26$ and $17$ actual ration boxes respectively. The proof of that is left as an exercise to the reader. For attaining the lower bound of $16\frac{23}{35}$ days, it is important that the two ration boxes which Professor Walkingholme can carry are not replaced by a single ration box of double capacity.

\section{Inequalities for lower bounds}

Consider the following problem. Suppose that Professor Walkingholme has reached as far as position $\gamma$ until now, but not farther. Suppose that he is on position $e_0$, and that before the current moment, he has unpacked $l$ rations on positive positions, say on positions $e_1, e_2, \ldots, e_l$, such that
$$
0 < e_l \le e_{l-1} \le \cdots \le e_2 \le e_1
$$
Define
$$
e_{l+1} = e_{l+2} = \cdots = 0
$$
Suppose that on the current moment, there are rations on positive positions $f_1, f_2, \ldots, f_{l'}$ which Professor Walkingholme has dropped before the current moment, in addition to the used rations on positions $e_1, e_2, \ldots, e_l$. Define
$$
f_{l'+1} = f_{l'+2} = \cdots = 0
$$
and assume that 
$$
f_1 = \max\{f_1,f_2,f_3,\ldots\}
$$
Take $k \in \N$ and $j \in \N$ such that $k + j \ge 1$.

\begin{lemma} \label{far_e1_A}
Suppose that either $j = 0$ or $e_1 \ge f_1$. Let $t$ be the time taken until now in days. Then
$$
10 t \ge \gamma + \sum_{i=2}^k e_i + \sum_{i=1}^j f_j - \frac12 e_0
$$
\end{lemma}

\begin{proof}
Professor Walkingholme carried ration boxes to $e_1, e_2, \ldots, e_k$ and $f_1, f_2, \allowbreak \ldots, f_j$, as far as these positions are nonzero. In addition, he walked from $e_1$ to $\gamma$. This adds up to
$$
\gamma + \sum_{i=2}^k e_i + \sum_{i=1}^j f_j
$$
miles in forward direction. The amount of walked miles in backward direction is $e_0$ miles less than that in forward direction. Since the total amount of walked miles is $20t$, the inequality follows.
\end{proof}

\begin{lemma} \label{far_f1_A}
Suppose that either $k = 0$ or $f_1 \ge e_1$. Let $t$ be the time taken until now in days. Then
$$
10 t \ge \gamma + \sum_{i=1}^k e_i + \sum_{i=2}^j f_j - \frac12 e_0
$$
\end{lemma}

\begin{proof}
The proof is similar to that of lemma \ref{far_e1_A}.
\end{proof}

Let $\beta = e_{2k+j} + 10$. During the 10 miles before eating a banana on position $e_i$, Professor Walkingholme can maximize the portion of his walk that is done in the interval $[\beta,\infty)$ by walking from position $e_i + 10$ to position $e_i$. He walks at most $\max\{e_i+10,\beta\} - \max\{e_i,\beta\}$ miles in $[\beta,\infty)$. During the 10 miles after eating a banana on position $e_i$, Professor Walkingholme can maximize the portion of his walk that is done in the interval $[\beta,\infty)$ by walking from position $e_i$ to position $e_i + 10$. Again, he walks at most $\max\{e_i+10,\beta\} - \max\{e_i,\beta\}$ miles in $[\beta,\infty)$.

Let $f_0 = \max\{e_0+10,\beta\}$. For $1 \le i \le k$, we estimate $\max\{e_i+10,\beta\} - \max\{e_i,\beta\}$ by $10$. For $k+1 \le i \le 2k+j$, we estimate $\max\{e_i+10,\beta\} - \max\{e_i,\beta\}$ by $\max\{e_i+10,\beta\} - \beta = e_i + 10 - \beta$. So Professor Walkingholme has walked 
$$
\alpha \le \big(f_0 - \max\{e_0,\beta\}\big) + 2 \sum_{i=1}^k 10 + 2 \sum_{i=k+1}^{2k+j} (e_i + 10 - \beta)
$$
miles in $[\beta,\infty)$. These $\alpha$ miles add up to a walk from position $\beta$ to position $\max\{e_0,\beta\}$, of which 
$$
\frac{\alpha - (\max\{e_0,\beta\} - \beta)}{2}
$$
miles are in backward direction and
\begin{align*}
\alpha' &= \frac{\alpha + (\max\{e_0,\beta\} - \beta)}{2} \\
&\le \frac12\big(f_0 - \beta\big) + \sum_{i=1}^k 10 + \sum_{i=k+1}^{2k+j} (e_i + 10 - \beta) \\
&= \frac12\big(f_0 - 10\big) + \sum_{i=k+1}^{2k+j-1} e_i + \frac12 e_{2k+j} + (2k+j) \cdot 10 - (k+j) \cdot \beta
\end{align*}
miles are in forward direction.

\begin{lemma} \label{conditions}
If either $e_0 \ge e_{2k+j}$, or $j \ge 1$ and $e_0 + 10 \ge f_j$, then
$$
\alpha' \le \frac12 e_0 + \sum_{i=k+1}^{2k+j-1} e_i + \frac12 e_{2k+j} + (2k+j) \cdot 10 - (k + j - 1) \cdot \beta - \min\{f_j,\beta\}
$$
\end{lemma}

\begin{proof}
If $e_0 \ge e_{2k+j}$, then $f_0 = e_0 + 10$, and the claim follows. So assume that $j \ge 1$ and $e_0 + 10 \ge f_j$. Then
\begin{align*}
\frac12(f_0 - 10) &= \frac12\big(e_0 + 10 + \beta - \min\{e_0+10,\beta\} - 10\big) \\
&\le \frac12\big(e_0 + \beta - \min\{f_j,\beta\}\big) \\
&\le \frac12 e_0 + \beta - \min\{f_j,\beta\}
\end{align*}
and the claim follows as well.
\end{proof}

\begin{lemma} \label{far_e1_B}
Suppose that the conditions of lemma \ref{conditions} hold. Suppose in addition that $k \ge 1$, and either $j = 0$ or $e_1 \ge f_1$. Then
$$
\gamma + \sum_{i=2}^k e_i + \sum_{i=1}^j f_i \le \frac12 e_0 + \sum_{i=k+1}^{2k+j-1} e_i + \frac12 e_{2k+j} + (2k+j) \cdot 10
$$
\end{lemma}

\begin{proof}
Professor Walkingholme carried ration boxes to positions $e_1, e_2, \ldots, e_k$ and $f_1, f_2, \allowbreak \ldots, f_k$, and also walked from position $e_1$ to position $\gamma$. If $j \ge 1$, then
\begin{align*}
\alpha' &\ge \big(\gamma - \min\{\gamma,\beta\}\big) + \sum_{i=2}^k \big(e_i - \min\{e_i,\beta\}\big) + \sum_{i=1}^j \big(f_i - \min\{f_i,\beta\}\big) \\
&\ge \gamma + \sum_{i=2}^k e_i + \sum_{i=1}^j f_i - (k + j - 1) \cdot \beta - \min\{f_j,\beta\}\big)
\end{align*}
If $j = 0$, then the above holds as well, because $\min\{f_0,\beta\} = \beta$. The claimed inequality follows by combining the above lower bound of $\alpha'$ with the upper bound in lemma \ref{conditions}.
\end{proof}

\begin{lemma} \label{far_f1_B}
Suppose that the conditions of lemma \ref{conditions} hold, and that $\gamma \le \frac12 e_0 + \frac12 e_1 + 10$. Suppose in addition that $j \ge 1$, and either $k = 0$ or $f_1 \ge e_1$. Then
$$
\gamma + \sum_{i=1}^k e_i + \sum_{i=2}^j f_i \le \frac12 e_0 + \sum_{i=k+1}^{2k+j-1} e_i + \frac12 e_{2k+j} + (2k+j) \cdot 10
$$
\end{lemma} 

\begin{proof}
If $j = 1$ and $k = 0$, then the claimed inequality comes down to $\gamma \le \frac12 e_0 + \frac12 e_1 + 10$, which was assumed.

Professor Walkingholme carried ration boxes to $e_1, e_2, \ldots, e_k$ and $f_1, f_2, \ldots, \allowbreak f_k$, and also walked from $f_1$ to $\gamma$. If $j \ge 2$, then
\begin{align*}
\alpha' &\ge \big(\gamma - \min\{\gamma,\beta\}\big) + \sum_{i=1}^k \big(e_i - \min\{e_i,\beta\}\big) + \sum_{i=2}^j \big(f_i - \min\{f_i,\beta\}\big) \\
&\ge \gamma + \sum_{i=1}^k e_i + \sum_{i=2}^j f_i + (k + j - 1) \cdot \beta - \min\{f_j,\beta\}
\end{align*}
If $j = 1$ and $k \ge 1$, then the above holds as well, because $f_1 \ge e_1$ and
$$
\alpha' \ge \gamma + \sum_{i=1}^k e_i + \sum_{i=2}^j f_i + (k + j - 1) \cdot \beta - \min\{e_1,\beta\}
$$
The claimed inequality follows by combining the above lower bound of $\alpha'$ with the upper bound in lemma \ref{conditions}.
\end{proof}

The above lower bounds were inspired by lower bounds in the solution for one way trips which take a whole number of days in \cite{cb}. The solution in \cite{cb} can easily be extended to one way trips which do not take a whole number of days. But as far as I know, there is no general solution for round trips. The above lower bounds could be of value for a possible solution. In \cite{cbl}, the same problem is solved as in \cite{cb}.

In the 2010 version of the paper, $\max\{e_0+10,\beta\} - \max\{e_0,\beta\}$ is estimated incorrectly. Fixing this led to the conditions in lemma \ref{conditions}. It is clear that the conditions in lemma \ref{conditions} may be assumed if the estimates are used for one way trips. The following lemma shows that the conditions of lemma \ref{conditions} may be assumed as well if the estimates are used for round trips. This lemma will not be used in this paper.

\begin{lemma}
Suppose that $e_0 \ge f_1$, and that in the last $\gamma - e_0$ miles, Professor Walkingholme walked from position $\gamma$ to position $e_0$ without unsealing a ration box between position $e_0$ and position $\gamma$. Then we may assume that the conditions of lemma \ref{conditions} are satisfied.
\end{lemma}

\begin{proof}
If $j \ge 1$, then $e_0 + 10 > f_1 \ge f_j$. So assume that $j = 0$. Then $k \ge 1$ and $e_2 \ge e_{2k+j}$. So it suffices to show that we may assume that $e_0 \ge e_2$.

We modify the strategy of Professor Walkingholme up to the current moment as follows. The last $2(\gamma-e_0)$ miles become a walk from position $e_0$ to position $\gamma$, followed by the walk from position $\gamma$ to position $e_0$ that Professor Walkinghole already did. Before the last $2(\gamma-e_0)$ miles, Professor Walkingholme must not go farther than $e_0$. So instead of being on position $v$, Professor Walkingholme is on position $\min\{v,e_0\}$.

The last $2(\gamma-e_0)$ miles were already on positions $\ge e_0$, so the transportations on positions $\le e_0$ are not affected. Since $\gamma-e_0 \le 20$, Professor Walkingholme unseals at most one ration box in $(e_0,\gamma]$ in the first half of the last $2(\gamma-e_0)$ miles. By assumption, he does not unseal a ration box in $(e_0,\gamma)$ in the second half of the last $2(\gamma-e_0)$ miles. So he needs to take at most two ration boxes with him in the last $2(\gamma-e_0)$ miles. Therefore, Professor Walkingholme can walk the last $2(\gamma-e_0)$ miles as indicated, and we may assume that $e_0 \ge e_2$.
\end{proof}

\section{Optimality of algorithms}

We computed the lower bounds as follows. We turned the used inequalities into equalities by adding positive symbolic constants to the smaller sides. We solved these inequalities with Maple 8. We used a Python script to construct algorithms in which the lower bounds are attained.

Let us first consider the method of making the circuit in the same way as crossing a distance of $100$ miles. We consider the last moment of this crossing: Professor Walkingholme is on position $e_0$ and there are no ration boxes left. We used the inequalities
\begin{align*}
\gamma &\ge e_0 & e_0 &\le e_1 + 20
\end{align*}
Furthermore, we used lemma \ref{far_e1_B} for $k = 1,2,\ldots,41$ and $j = 0$, and lemma \ref{far_e1_A} for $k = 42,43,\ldots,82$ and $j = 0$. We obtained that
\begin{equation} \label{est_O}
t_O \ge \frac{1536}{275}\gamma - \frac{2382}{5}
\end{equation}
Substituting $\gamma = 100$ in \eqref{est_O} yields $t \ge 82\frac{8}{55}$. This bound can be attained if we allow Professor Walkingholme to start on a specific moment during the day, to count the first day partially. 

Let us next consider the method of making the circuit by way of two round trips. Let $d$ be a day on which Professor Walkingholme reaches $\gamma$ in one of both round trips. We consider the moment at the end of day $d$: Professor Walkingholme is on position $e_0'$ and there are ration boxes on positions $f_1, f_2, \cdots, f_{l'}$. It is possible that $e_0' < 0$, i.e.\@ $e_0'$ is on the other side of the starting point. But walking on negative positions does not have any purpose for the round trip at hand. So we may assume that $e_0' \ge 0$. 

There are two cases: $e_1 \ge f_1$ and $f_1 \ge e_1$. For the case $e_1 \ge f_1$, we used the equality $f_1 = e_0'$ and the inequalities
\begin{align*}
f_2 &\ge f_1 - 20 & \gamma &\le \frac12 e_0' + \frac12 e_1 + 10 & e_1 \le e_2 + 20
\end{align*}
Furthermore, we used lemma \ref{far_e1_B} for $e_0 = e_0'$, $k = 2,3,4$ and $j = 1$, and for $e_0 = e_0'$, $k = 4,5,\ldots,8$ and $j = 2$. In addition, we used lemma \ref{far_e1_A} for $e_0 = 0$, $k = 9,10,\ldots,18$ and $j = 2$. We obtained that
\begin{equation} \label{est_R1}
t_R \ge \frac{27}{20} \gamma - \frac{375}8
\end{equation}
Using lemma \ref{far_e1_A} with $e_0 = 0$ is wrong in some sense, because if $e_0' > 0$, then Professor Walkingholme will be on position $e_0 = 0$ on a later moment. In the meantime, the ration boxes on positions $f_1, f_2, \ldots, f_j$ may be unsealed. But this does not affect the formula of lemma \ref{far_e1_A}, provided the ration boxes on positions $f_1, f_2, \ldots, f_j$ are not moved any more before unsealing. One may indeed assume that the ration boxes are not moved any more, just as ration boxes are not moved in part C in the proof of lemma \ref{ABC}. 

Substituting $\gamma = 50$ in \eqref{est_R1} yields $t \ge 20\frac{5}{8}$. This bound can be attained even if Professor Walkingholme must start at dawn. But since the round trip does not take a whole number of days, it cannot be followed by a second round trip of the same type.

For the case $f_1 \ge e_1$, we used the equality $f_1 = e_0'$ and the inequalities
\begin{align*}
f_2 &\ge f_1 - 20 & f_3 &\ge f_2 - 20 & f_1 &\le \gamma
\end{align*}
We computed two estimates of $t_R$. For the first estimate, we used lemma \ref{far_f1_B} for $e_0 = e_0'$, $k = 0$ and $j = 1$, and for $e_0 = e_0'$, $k = 1,2,\ldots,7$ and $j = 2$. In addition, we used lemma \ref{far_e1_A} for $e_0 = 0$, $k = 8,9,\ldots,16$ and $j = 2$. We obtained that
\begin{equation} \label{est_R2}
t_R \ge \frac{181}{140} \gamma - 44
\end{equation}
For the second estimate, we used lemma \ref{far_f1_B} for $e_0 = e_0'$, $k = 0$ and $j = 1$, for $e_0 = e_0'$, $k = 1,2,3$ and $j = 2$, and for $e_0 = e_0'$, $k = 4,5,\ldots,7$ and $j = 3$. In addition, we used lemma \ref{far_e1_A} for $e_0 = 0$, $k = 8,9,\ldots,17$ and $j = 3$. We obtained that
\begin{equation} \label{est_R3}
t_R \ge \frac{191}{140} \gamma - \frac{333}7
\end{equation}
Substituting $\gamma = 50$ in both \eqref{est_R2} and \eqref{est_R3} yields $t \ge 20\frac{9}{14}$ which indicates that $e_1 \ge f_1$ is the better option for a round trip of $50$ miles. But the round trips do not need to be $50$ miles: they must only add up to $100$ miles. So additional argumentation is needed. Let
$$
\gamma' = 100 - \gamma
$$
From \eqref{est_R2} and \eqref{est_R3}, it follows that
\begin{equation} \label{est_R23}
t_R \ge \frac15 \Big(\frac{181}{140} \gamma - 44\Big) + \frac45 \Big(\frac{191}{140} \gamma - \frac{333}7\Big) = \frac{27}{20} \gamma - \frac{328}{7}
\end{equation}
So if $f_1 \ge e_1$, then a round trip to $\gamma$ takes at least $-\frac{328}{7} - -\frac{375}{8} = \frac{1}{56}$ day more than the right hand side of $\eqref{est_R1}$. Therefore, at least
$$
\Big(\frac{27}{20} \gamma - \frac{375}8\Big) + \Big(\frac{27}{20} \gamma' - \frac{375}8\Big) = 41\frac14
$$
days are necessary for both round trips. 

The above lower bound can be attained if Professor Walkingholme walks $13\frac13$ miles on the first day, so the first day counts for $\frac23$ days. On the first day, he carries two ration boxes $3\frac13$ miles from the starting point, but on different sides, and returns to the starting point. On day $2, 3, \ldots, 22$, he makes a round trip of $50\frac{5}{18}$ miles, returning with $\frac13$ ration box. He uses the $\frac13$ ration box to carry another ration box $3\frac13$ miles from the starting point, and to return to the starting point afterwards. On day $23, 24, \ldots, 42$, he makes a round trip of $49\frac{13}{18}$ miles, returning with $\frac{5}{12}$ ration box. So on day $1$ and day $42$ together, Professor Walkingholme uses $\frac23 + \frac{7}{12} = 1\frac14$ ration boxes. 

Let us finally consider the method of using part A, part B and part C as indicated in lemma \ref{ABC}. We computed two estimates for part A. Let $e_1, e_2, e_3, \ldots$ be the positions where a ration box is unsealed before reaching position $\gamma'$ for the last time in part A or on the moment of reaching position $\gamma'$ for the last time in part A. Let $f_2, f_4, f_6, \ldots$ be the positions where a ration box is unsealed after reaching position $\gamma'$ for the last time in part A. Let $f_1, f_3, f_5, \ldots$ be the positions where the ration boxes which are unsealed in part C are at the end of part A. 

Let $d$ be the last day on which Professor Walkingholme is on position $\gamma'$ in part
A. We consider the moment at the end of day $d$: we may assume that Professor Walkingholme is on a position $e_0' \ge 0$ and that there are ration boxes on positions $f_1$ and $f_2$, because of the following. The ration box on position $\gamma'$ will not be moved any more before unsealing. If $e_0' > 0$, then a ration box needs to be unsealed on position $e_0'$ the next morning. So we may assume that $f_1 = \gamma'$ and $f_2 = e_0'$.

For the first estimate for part A, we used the (in)equalities
\begin{align*}
f_1 &= \gamma' & f_2 &= e_0' & f_3 &\ge f_1 - 20 & f_4 &\ge f_2 - 20 & \gamma' &\le \frac12 e_0' + \frac12 e_1 + 10
\end{align*}
Furthermore, we used lemma \ref{far_f1_B} for $e_0 = e_0'$, $k = 0,1$ and $j = 2$. In addition, we used lemma \ref{far_f1_A} for $e_0 = 0$, $k = 2,3$ and $j = 3$, and for $e_0 = 0$, $k=4$ and $j=4$. We obtained that
\begin{equation} \label{est_A1}
t_A \ge \frac{7}{10} \gamma' - 11
\end{equation}
Substituting $\gamma = 70$ in \eqref{est_A1} yields $t \ge 10$, so part A of Algorithm \ref{dudeney} is optimal. Substituting $\gamma = 71\frac14$ in \eqref{est_A1} yields $t \ge 9\frac18$, so part A of Algorithm \ref{best} is optimal as well.

For the second estimate for part A, we used the (in)equalities
\begin{align*}
f_1 &= \gamma' & f_2 &= e_0' & f_3 &\ge f_1 - 20 & \gamma' &\le \frac12 e_0' + \frac12 e_1 + 10
\end{align*}
Furthermore, we used lemma \ref{far_f1_B} for $e_0 = e_0'$, $k = 0,1$ and $j = 2$. In addition, we used lemma \ref{far_f1_A} for $e_0 = 0$, $k = 2,3,4$ and $j = 3$. We obtained that
\begin{equation} \label{est_A3}
t_A \ge \frac{22}{35} \gamma' - \frac{64}{7}
\end{equation}
Substituting $\gamma = 75\frac{20}{29}$ in \eqref{est_A3} yields $t \ge 6\frac{4}{29}$, so part A of Algorithm \ref{dawn} is optimal.

We computed three estimates for part B. We consider the last moment of part B: Professor Walkingholme is on position $e_0$ and there are no ration boxes from the side of part B left. We used the inequalities
\begin{align*}
\gamma &\ge e_0 & e_0 &\le e_1 + 20
\end{align*}
Furthermore, we used lemma \ref{far_e1_B} for $k = 1,2,\ldots,n$ and $j = 0$, and lemma \ref{far_e1_A} for $k = n+1,n+2,\ldots,2n$ and $j = 0$. 
\newcounter{foo}%
\setcounter{foo}{\value{equation}}%
\addtocounter{foo}{-4}%
We obtained inequality ($n + \thefoo$) of 
\begin{align}
t_B &\ge \frac{24}{35} \gamma - \frac{258}7 \label{est_B1} \\
t_B &\ge \frac{4}{5} \gamma - 45 \label{est_B2} \\
t_B &\ge \frac{24}{25} \gamma - \frac{284}5 \label{est_B3} 
\end{align}
Substituting $\gamma = 70$ in \eqref{est_B1} yields $t_B \ge 11\frac{1}{7}$. This bound can be attained if we allow Professor Walkingholme to start on a specific moment during the day, to count the first day partially. So part B of Algorithm 1 is not optimal. Substituting $\gamma = 71\frac14$ in both \eqref{est_B1} and \eqref{est_B2} yields $t_B \ge 12$, so part B of Algorithm 2 is optimal. Substituting $\gamma = 75\frac{20}{29}$ in \eqref{est_B3} yields $t_B \ge 15\frac{25}{29}$, so part B of Algorithm 3 is optimal.

We next prove the optimality of Algorithm 2. For part C, we have
\begin{equation} \label{est_C}
t_C \ge \frac{1}{20} \gamma'
\end{equation}
From \eqref{est_B1} and \eqref{est_B2}, it follows that
\begin{equation} \label{est_B12}
t_B \ge \frac{7}{16} \Big(\frac{24}{35} \gamma - \frac{258}7\Big) 
       + \frac{9}{16} \Big(\frac{4}{5} \gamma - 45\Big) 
       = \frac{3}{4} \gamma - \frac{663}{16}
\end{equation}
From \eqref{est_A1}, \eqref{est_B12} and \eqref{est_C}, it follows that
$$
t_A + t_B + t_C \ge \frac{361}{16} = 22\frac{9}{16}
$$
which proves the optimality of Algorithm 2.

Finally, we prove the optimality of Algorithm 3. From \eqref{est_A1}, \eqref{est_B1} and \eqref{est_B2}, it follows that
\begin{equation} \label{est_B21}
t_A + t_B \ge \Big(\frac{7}{10} \gamma' - 11\Big) + \frac{7}{8} \Big(\frac{24}{35} \gamma - \frac{258}7\Big) 
       + \frac{1}{8} \Big(\frac{4}{5} \gamma - 45\Big) 
       = 21\frac{1}{8}
\end{equation}
If part A and part B must take a whole amount of days, then they take at least $22$ days. From \eqref{est_A3}, \eqref{est_B3} and \eqref{est_C}, it follows that
$$
t_A + t_B + t_C \ge \frac{197}{700}\gamma + \frac{67}{35} 
= \frac{197}{700} \Big(\gamma - 75\frac{20}{29}\Big) + 23\frac{25}{116}
$$
Since part A and part B together take $22$ days in Algorithm 3 and at least $22$ days in general, Professor Walkingholme can only improve on Algorithm 3 with a larger value of $\gamma$. So we may assume that $\gamma \ge 75\frac{20}{29}$. Therefore $t_A + t_B + t_C \ge 23\frac{25}{116}$, so Algorithm 3 is optimal.

\end{document}